\documentclass[notitlepage,leqno,11pt]{article}
\usepackage{amssymb}
\catcode`\@=11 \@addtoreset{equation}{section}

\catcode`\@=12

\usepackage{latexsym}
\usepackage{textcomp}
\usepackage{amsmath}
\usepackage{amsfonts}
\usepackage{amssymb}
\usepackage{mathrsfs}

\renewcommand{\a }{\alpha }
\renewcommand{\b }{\beta }
\renewcommand{\d}{\delta }

\newcommand{\D }{\Delta }

\newcommand{\e }{\varepsilon }
\newcommand{\g }{\gamma}

\newcommand{\G }{\Gamma}
\renewcommand{\l }{\lambda }

\newcommand{\n }{\nabla }
\newcommand{\vp }{\varphi }

\newcommand{\s }{\sigma }
\newcommand{\Sig }{\Sigma}

\renewcommand{\th }{\theta }
\renewcommand{\o }{\omega }
\renewcommand{\O }{\Omega }

\newcommand{\ov}{\overline}

\newcommand{\be}{\begin{equation}}
\newcommand{\ee}{\end{equation}}

\newcommand{\R}{\mathbb{R}}
\newcommand{\N}{\mathbb{N}}
\renewcommand{\S}{\mathbb{S}}

\newcommand{\C}{\mathbb{C}}

\newcommand{\de}{\partial}

\newcommand{\ti}{\widetilde}

\newcommand{\M}{\mathcal{M}}

\newcommand{\ra}{{\rangle}}
\newcommand{\la}{{\langle}}

\newcommand{\calO }{\mathcal{O}}
\newcommand{\calC }{\mathcal{C}}

\newcommand{\calD }{\mathcal{D}}

\newcommand{\calB }{\mathcal{B}}

\newcommand{\calU}{{\mathcal U}}
\newcommand{\calM}{{\mathcal M}}

\newtheorem{Theorem}{Theorem}[section]
\newtheorem{Lemma}[Theorem]{Lemma}
\newtheorem{Proposition}[Theorem]{Proposition}
\newtheorem{Corollary}[Theorem]{Corollary}
\newtheorem{Remark}[Theorem]{Remark}

\def\proof{\noindent{{\bf Proof. }}}
\def\square{\vbox{
    \hrule height .4pt
    \hbox{\vrule width .4pt height 7pt \kern 7pt
       \vrule width .4pt}
    \hrule height .4pt }}

\def\square{\vbox{
    \hrule height .4pt
    \hbox{\vrule width .4pt height 7pt \kern 7pt
       \vrule width .4pt}
    \hrule height .4pt }}

\def\QED{\hfill {$\square$}\goodbreak \medskip}

\def\R{{\mathbb R}}

\def\C{{\mathcal C}}
\def\S{{\mathbb S}}
\def\f{{\varphi}}

\def\eps{{\varepsilon}}

\font\sc=cmcsc9  \textwidth=14truecm
\begin{document}

\title{ Nonexistence of distributional supersolutions of a semilinear elliptic equation with
Hardy potential}

\author{Mouhamed Moustapha Fall\footnote
{\footnotesize{Institut f\"{u}r Mathematik Goethe-Universit\"{a}t
Frankfurt, Robert-Mayer-Str. 10, D-60054 Frankfurt am Main, Germany.
 E-mail: {\tt fall@math.uni-frankfurt.de,  mouhamed.m.fall@gmail.com}. } } }
\date{}

\maketitle

\bigskip

\noindent {\footnotesize{\bf Abstract.} In this paper we study nonexistence of
non-negative distributional supersolutions for a class of
semilinear elliptic equations involving  inverse-square potentials.}
\bigskip\bigskip

\noindent{\footnotesize{{\it Key Words:} Hardy inequality, critical
exponent, nonexistence, distributional solutions, Fermi coordinates, Emden-Flower transform.}}

\noindent{\footnotesize{{\it 2000 Mathematics Subject
Classification:} 35D05, 35J15.}}

\section*{Introduction}
Let $\O$ define a domain of  $\R^N$, $N\ge 3$. In this paper, we study
nonnegative functions $u $ satisfying
\begin{equation}
\label{eq:inequalityb}
-\Delta u-b(x)\, u \ge u^p\quad\textrm{ in } \calD'(\O),
\end{equation}
with $p>1$, $b \gneqq0$ and $b \in L^1_{loc}(\O) $ is a singular potential of Hardy-type.
More precisely, we are interested in distributional solutions to (\ref{eq:inequalityb}), that is,
functions $u\in L^p_{\rm loc}(\O)$ such that $b(x)\,u\in L^1_{\rm loc}(\O)$ and
$$
\int_\O u(-\Delta\f-  b(x)\f)~dx\ge\int_\O u^p\f~dx
\quad
\textrm{ $\forall\f\in C^\infty_c(\O)$, $\f\ge 0$.}
$$
 The study of nonexistence results  of (very) weak solution to problem \eqref{eq:inequalityb}  goes back to \cite{BC}, where the authors
 were motivated by the failure of the Implicit Function Theorem. Further
 references in this direction are \cite{BDT}, \cite{DD-H}, \cite{D}, \cite{DN}. We also quote \cite{AS-CPDE}, \cite{AS-RIMS} \cite{PoTe}, \cite{T}, 
\cite{LLM-Proce}, \cite{KLS},
\cite{KLS-JDE}, \cite{KLS-Trans}.\\
In this paper, we study nonexistence of solutions to \eqref{eq:inequalityb}
 when $\de\O$ possesses a conical  singularity at $0$ as well as when $\de\O$
is of class $ C^2$ at 0. Higher dimensional singularity will be also considered.\\
For any   domain $\Sigma$
in the unit sphere $\S^{N-1}$ we introduce the cone
$$
\C_\Sigma:=\left\{ r\sigma\in\R^N~|~r>0~,~\sigma\in\Sigma~\right\}~\!.
$$
We recall that  the best constant in the Hardy inequality for functions supported by
$\C_\Sigma$ is given by
$$
\mu(\C_\Sigma):= \inf_{u\in C^\infty_c(\C_\Sigma)}
\frac{\displaystyle\int_{\C_\Sigma}|\n u|^2~dx}
 {\displaystyle\int_{\C_\Sigma}|x|^{-2}u^2~dx}=
 \frac{(N-2)^2}{4}+\lambda_1(\Sigma)~\!,
 $$
 where $\lambda_1(\Sigma)$ is the first Dirichlet eigenvalue for the
 Laplace-Beltrami operator on $\Sigma$  (\cite{FaMu1}, \cite{PT}).
 For a given radius $R>0$ we introduce the
 cone-like domain
 $$
 \C_\Sigma^R:=\C_\Sigma\cap B_R=\left\{ r\sigma~|~r\in(0,R)~,~\sigma\in\Sigma~\right\},
 $$
 where $B_R$ is the ball of radius $R$
centered at $0$.
We study the inequality
\begin{equation}
\label{eq:inequalityC}
-\Delta u-\frac{c}{|x|^2}\, u \ge u^p\quad\textrm{ in } \calD'\left(\C_\Sigma^R \right),
\end{equation}
with
$$
 \lambda_1(\Sigma)<c\le \mu(\C_\Sigma)~\!.
$$
By homogeneity, an important   role is played by
$$
\alpha^-_\Sigma:=\frac{N-2}{2}-\sqrt{\mu( \calC_\Sig)-c}~\!,
$$
which is the smallest root of the equation
$$
\a^2-(N-2)\a+c-\l_1(\Sig)=0.
$$
We notice that the restriction $c\leq \mu(\calC_\Sig)$ is not restrictive (see Remark \ref{re:ne-cgmu} below) and in addition
 $\alpha^-_\Sigma>0$ when $ c>  \lambda_1(\Sigma)$. Finally we define
$$p_\Sig=1+\frac{2}{\alpha^-_\Sigma}~\!.$$
We observe that
$p_\Sig=\frac{N+2}{N-2}$ when $c=\mu(\calC_\Sig)$ while $p_\Sig>\frac{N+2}{N-2} $ as soon as $c<\mu(\calC_\Sig)$.\\
In  \cite{BDT}, the authors
have studied the case $\O=B_R\setminus\{0\}= \C_{\S^{N-1}}^R$.
They   proved  that \eqref{eq:inequalityC} has a non-trivial solution
in $B_R\setminus\{0\}$ if and only if $p< p_{\S^{N-1}}$.\\
Our first result generalizes   the nonexistence result  in \cite{BDT}  to  cone-like domains.
\begin{Theorem}\label{th:neC}
Let $\calC_{\Sig}^R$ be a cone-like domain of $\R^N,$ $N\geq3$. For $\lambda_1(\Sigma)<c\le \mu(\C_\Sigma)$, let $u\in L^p_{loc}\left(\calC_\Sig^R\right)$
be  non-negative such that
$$
-\D u-\frac{c}{|x|^2}u\geq u^p\quad\textrm{ in }\,
\calD'\left(\calC_\Sig^R\right)~\!.
$$
If $p\geq p_\Sig$  then $u\equiv0$.
\end{Theorem}

Theorem \ref{th:neC} improves  a part of  the nonexistence
results obtained in \cite{LLM}, where more regular supersolutions were considered.
We notice that the assumption $p\ge p_\Sigma$ is sharp
 (see the existence result in \cite{LLM}, Theorem 1.2). \\
\medskip
We next consider the case where $0\in \de \O$  with $\de\O$ is smooth at   $0$ and $b(x)=c|x|^{-2}$. We define
$$
\mu(\O):= \inf_{u\in C^\infty_c(\O)}
\frac{\displaystyle\int_{\O}|\n u|^2~dx}
 {\displaystyle\int_{\O}|x|^{-2}u^2~dx}.
 $$
Put $\O_r:=\O\cap B_r(0)$. Recently it was proved in \cite{mmf} that, there exits $r_0=r_0(\O)>0$ such that for
all $r\in(0,r_0)$
\begin{equation}
\label{eq:muOr}
\mu(\O_r)=\mu\left(\calC_{\S^{N-1}_+}\right)=\frac{N^2}{4}~\!,
\end{equation}
with $\S^{N-1}_+$ is a hemisphere centered at 0 so that $ \calC_{\S^{N-1}_+}$ is a half-space. 
We have obtained:
\begin{Theorem}\label{th:ne-sm-smoo-loc} Let $\O$ be a smooth domain of $\R^{N}$, $N\geq3$,
with $0\in\de\O$. Let $r>0$ small so that \eqref{eq:muOr} holds.
  For   $N-1<c\leq \frac{N^2}{4}$,
let $u\in L^p_{loc}\left(\O_r\right)$
be  non-negative such that
$$
-\D u-\frac{c}{|x|^2}u\geq u^p\quad\textrm{ in }\,
\calD'\left(\O_r\right)~\!.
$$
If $p\geq p_{\S^{N-1}_+} $  then $u\equiv0$.
\end{Theorem}
Here also the nonexistence of nontrivial solution for $c\in(N-1,N^2/4]$ is sharp, see Proposition \ref{prop:locex}.\\
When we consider general domains, we face some obstacles in
the restriction of the parameter $c$. This is due to the fact that $\mu(\O)$  is not in general smaller than $N-1$ for  smooth domains $\O$,
with $0\in \de\O$,  see \cite{FaMu1}. 
A consequence of Theorem \ref{th:ne-sm-smoo-loc} is:
\begin{Corollary}\label{cor:ne-sm}
Let $\O$ be a smooth domain of $\R^{N}$, $N\geq3$,
with $0\in\de\O$.
Assume that $N-1< c \leq\mu(\O)$. Suppose that there exists
$u\in L^p_{loc}(\O)$, $u\geq0$ such that
$$
-\D u-\frac{c}{|x|^{2}}u\geq u^p\quad \textrm{ in }\calD'(\O).
$$
If   $p\geq p_{\S^{N-1}_+}$ then $u\equiv0$.
\end{Corollary}
In Corollary \ref{cor:ne-sm} above, we assume that the interval $[N-1,\mu(\O)]$ is not empty. This is not in general true  (see Remark \ref{rm:moss} below).
However it holds for various domains or in higher dimensions. Indeed,
%
we first observe that the inequality $\frac{(N-2)^2}{4}< \mu(\O)\leq\frac{N^2}{4}$ is
 valid for every smooth bounded domain $\O$ with $0\in \de\O$, see \cite{FaMu1}.
In particular $ \mu(\O)>N-1$  whenever $N\geq 7$. Hence  we get:
\begin{Corollary}\label{cor.N7}
 Let $\O$ be a smooth bounded domain of $\R^N$,  $N\geq 7$, with $0\in\de\O$.
Let $ N-1<c\leq \mu(\O)$
 and  $u\in L^p_{loc}\left(\O\right)$
be  non-negative such that
$$
-\D u-\frac{c}{|x|^2}u\geq u^p\quad\textrm{ in }\,
\calD'\left(\O\right)~\!.
$$
If $p\geq p_{\S^{N-1}_+}$  then $u\equiv0$.
\end{Corollary}
 When  $\O$ is a smooth domain (not necessarily bounded), with $0\in\de\O$, is contained in the half-space $\calC_{\S^{N-1}_+}$
 then obviously  $\mu(\O)=\frac{N^2}{4} $ by  \eqref{eq:muOr}. In particular,
thanks to {Theorem} \ref{th:ne-sm-smoo-loc},  the restriction $N\geq 7$ in Corollary \ref{cor.N7} and the boundedness of $\O$
can be removed. Indeed, we have:
\begin{Corollary}\label{cor.Hs}
 Let $\O$ be a smooth  domain  of the half-space $\calC_{\S^{N-1}_+} $,  $N\geq 3$, with $0\in\de\O$.
Let $ N-1<c\leq \frac{N^2}{4} $
 and $u\in L^p_{loc}\left(\O\right)$
be non-negative such that
$$
-\D u-\frac{c}{|x|^2}u\geq u^p\quad\textrm{ in }\,
\calD'\left(\O\right)~\!.
$$
If $p\geq p_{\S^{N-1}_+}$  then $u\equiv0$.
\end{Corollary}
%
\begin{Remark}\label{rm:moss}
According to our argument, the assumption  $N-1< \mu(\O)$ is crucial
 because it implies that $1< p_{\S^{N-1}_+}< \infty$ when $c> N-1$.
However it  is not valid for every smooth domain.
In fact, one can construct  a family of
 smooth bounded domains $\O^\e$,  for which $\mu(\O^\e)\leq \frac{(N-2)^2}{4}+\e$, for  $\e>0$ small,
 see \cite{FaMu1}, \cite{mmf-p}.
\end{Remark}
\begin{Remark}
The conclusion in theorems \ref{th:neC}, \ref{th:ne-sm-smoo-loc} still holds
when $u^p$ is replaced by $|x|^s u^{q}$ with $\l_1(\Sig)< c \leq \mu\left(\calC_{\Sig}\right)$.
In this case one has to replace $p_{\Sig}$ with $q_{\Sig}=1+\frac{2-s}{\a^-_{\Sig}}$.
\end{Remark}

We prove our nonexistence results via a linearization argument which were also used in \cite{LLM}. However when working with
weaker notion of solutions, further analysis are required. Our approach  is to obtain a quite sharp lower estimate on $u$ in such
a way that $u^{p-1}$ is somehow  proportional to $b(x)$ and to look the problem as a linear problem: $-\D u- b(x)u-u^{p-1}\,u\geq0$ in $\calD'(\O)$.
This leads to the inequality (see Lemma \ref{lem:AP})
 \be\label{eq:apineq-i}
 \int_{\O}|\n\vp|^2-\int_{\O}b(x)\vp^2\geq \int_{\O}u^{p-1}\vp^2\quad\forall \vp\in C^\infty_c(\O).
 \ee
 By using appropriate test functions in \eqref{eq:apineq-i}, we were able to contradict the existence of solutions.
 To lower estimate $u$, we construct sub-solutions
for the operator $L:=-\D -b(x)$. On the other hand since we are working with "very weak" supersolutions in non-smooth domains,
and  the operator $L$ does not in general satisfies the maximum principle, we have proved a comparison principle (see Lemma \ref{lem:min-sol}
 in Section \ref{S:preliminaries}). We  achieve this by  requiring $L$ to be coercive. Since in this paper the potential $b(x)$ is of Hardy-type,
such coercivity is nothing but improvements of Hardy inequalities. The comparison principle  allows us to put below $u$ a  more regular function $v$.
Such function $v$ turns out to be a supersolution for $L$ and therefore  can be  lower estimated by the sub-solutions via standard arguments.\\

The paper is organized as follows. In Section \ref{S:preliminaries} we prove
some preliminary results, which are mainly  used in the paper.
The proofs of Theorems \ref{th:neC}, \ref{th:ne-sm-smoo-loc}
will be carried out in Sections \ref{S:cone},  \ref{S:smooth} respectively. Finally in the last section, we study the problem
\be\label{eq:petplmi}
\begin{cases}
-\D u - \frac{(N-k-2)^2}{4}\,\frac{1}{\textrm{dist}(x, \Gamma)^{2}}\, {q(x)}\,u\geq u^p\quad\textrm{ in } \calD'(\O\setminus \Gamma), \\
u\in L^p_{loc}(\O\setminus \Gamma),\\
u\gneqq0,
\end{cases}
\ee
where $\Gamma$ is a smooth closed submanifold of $\O$ and $q$ is a nonnegative weight.

\section{Preliminaries and comparison lemmata}
\label{S:preliminaries}
Let $\O$ be a  bounded open subset of $\R^N$. In this section we deal with comparison results
involving a differential operator of the type
$$
-\D-b(x)~\!,
$$
where $b\in L^1_{loc}(\O)$ is a given non-negative weight. We shall always assume that
$-\D-b(x)$ is coercive,  in the sense that there exists a
 constant $C(\O)>0$ such that
\begin{equation}\label{eq:coerciv}
 \int_{\O}|\n u|^2~dx-\int_{\O}b(x)u^2~dx\geq C(\O) \int_{\O}u^2~dx\quad\textrm{for any
 $u\in C^\infty_c(\O)$.}
\end{equation}
Following \cite{DD}, we define
the space  $H(\O)$ as the completion of  $C^\infty_c(\O)$ with respect to the
 scalar product
 $$
(u,v)\mapsto \int_{\O}\n u \n v~ dx-\int_{\O}b(x)uv~ dx~\!.
 $$
The scalar product in $H(\O)$ will be denoted by $\la \cdot ,\cdot \ra_{H(\O)}$.\\
Clearly $H^1_0(\O) \hookrightarrow H(\O)\hookrightarrow L^2(\O)$ by (\ref{eq:coerciv}), and hence
$L^2(\O)$ embeds into the dual space $H(\O)'$. 
By the
Lax Milligram theorem, for any  $f\in L^2(\O)$
there exists a unique function $v\in H(\O)$ such that
$$
-\Delta v - b(x) v=f\quad\textrm{in $H(\O)$,}
$$
that is,
$$
\la v,\vp\ra_{H(\O)}=\int_{\O} f \vp~dx\quad\textrm{for any
$\vp\in H(\O)$.}
$$

\begin{Remark}
\label{R:H2} Observe that if $b\in L^\infty(\O)$ then $H(\O)= H^1_0(\O)$ since
$$
C\int_\O|\nabla u|^2~dx\le\|u\|^2_{H(\O)}\le\int_\O|\nabla u|^2~dx~\!,
$$
where the constant $C>0$ depends only on $C(\O)$,
 and on the $L^\infty$ norm of $b$.
\end{Remark}

We start with the following technical result which will be useful in the sequel.
\begin{Lemma}\label{lem:min-solH10}
Let   $u\in
L^1_{loc}(\O)$ be non-negative and $g\in L^2(\O)$  such that
$$
-\D u\geq g \quad\textrm{ in } \calD'(\O).
$$
 Let $v\in
H^1_0(\O)$ be the solution to
$$
-\Delta v = g\quad\textrm{ in } \O.
$$
Then
 $$
v\leq u\quad\textit{ in $\O$.}
$$
\end{Lemma}
\proof For $\e>0$, define
$\O_\e=\{x\in\O\,:\,\textrm{dist}(x,\de\O)>\e\}$.
Let $\ti{\O}_\e$ be a smooth open set
 compactly contained in $\O$ and containing $\O_\e$.
Denote by $\rho_n$  the
standard mollifier and put $u_n=\rho_n*u$. Then for $\e>0$ there
exists $N_\e$ such that  $u_n$ is smooth in $\ti{\O}_\e$ up to the boundary for all $n\geq N_\e$.
Consider $v_{\e,n}\in H^1_0(\ti{\O}_\e)$ be the solution of
$ -\D v_{\e,n}=\rho_n* g=g_n $ in $\ti{\O}_\e.$
Clearly $-\D(u_n- v_{\e,n})\geq0$ in $\ti{\O}_\e$ and $u_n- v_{\e,n}\geq0$ on
$\de\ti{\O}_\e$, because $u$ is non-negative. It turns out that $u_n- v_{\e,n}\geq0$ in $\ti{\O}_\e$ by the maximum principle.
Letting $v_\e\in H^1_0(\ti{\O}_\e)$ be the solution of $-\D v_\e=g$ in $\ti{\O}_\e$, by H\"{o}lder and  Poincar\'e inequalities,
we have that $\|v_{\e,n}-v_\e\|_{ H^1_0(\ti{\O}_\e) }\leq C\|g_n-g\|_{L^2(\ti{\O}_\e)}$, with $C>0$ is a constant independent on $n$.
In particular  $v_{\e,n}$  converges to  $v_\e$  in $\ti{\O}_\e$. Therefore
 $u\ge v_\eps$  in $\ti{\O}_\e$. To conclude, it suffices to notice that
$v_\e \to v$ weakly in
$H^1_0(\O)$ and pointwise  in $\O$.
 \QED

We have the following comparison principle.
\begin{Lemma}\label{lem:min-sol}
Let $u\in L^1_{loc}(\O)$ be non-negative with $b(x)u\in L^1_{loc}(\O)$ and
let  $f\in L^2(\O)$ with $f\geq 0$ such that
$$
-\D u-b(x)u\geq f \quad\textit{ in $ \calD'(\O)$.}
$$
Let  $v\in H(\O)$ be the  solution of
$$
-\Delta v - b(x) v=f\quad\textrm{in $H(\O)$.}
$$
Then
 $$
v\leq u\quad\textit{  in $\O$.}
$$
\end{Lemma}
\proof
\textbf{Step 1:}\textit{ We first prove the result  if  $b\in L^\infty(\O)$.}\\
 We  let  $v_0\in H^1_0(\O)$ solving
$$
 -\D v_0= f\quad\textrm{ in $\O$}.
$$
Then $ 0\leq v_0\leq u$ in $\O$ by
Lemma \ref{lem:min-solH10} and because $f\geq 0$.
We define inductively the sequence $v_n\in H^1_0(\O)$ by
$$
 -\D v_1= {b}(x) v_{0}+f\quad\textrm{in $\O$},\quad\quad
 -\D v_n= {b}(x) v_{n-1}+f \quad\textrm{in $\O$}.
$$
Since ${b}\geq 0$, we have  $-\D u\geq  {b}(x)v_0+ f$ in $\calD'(\O) $. Thus   using  once again Lemma \ref{lem:min-solH10}, we  obtain
$ v_0\leq v_{1}\leq u $ in $\O$. By induction, we  have  
$$
v_0\leq v_1\leq \dots\leq  v_{n}\leq u\quad\textrm{  in $\O$}\quad \forall n\in\N.
$$
Since $v_{n-1}\leq v_n$ in $\O$, we have
$$
\int_{\O}|\n v_n|^2dx-\int_{\O}b(x)|v_n|^2\leq \int_{\O}f(x)v_ndx.
$$
By H\"{older} inequality and \eqref{eq:coerciv} (see Remark \ref{R:H2}) $v_n$ is bounded in $H^1_0(\O)$.
We conclude  that $v_n \rightharpoonup v$  in $H^1_0(\O)$ as $n\to \infty$ which is the unique solution to 
$$
-\D v= {b}(x) v+f\quad\textrm{in $\O$}.
$$
Since $v_n\to v$ in $L^2(\O)$, we get  $v\leq u$ in $\O$.

\bigskip 
\noindent
\textbf{Step 2:} \textit{Conclusion of the proof.}\\
We put $b_k(x)=\min(b(x),k)$ for every $k\in\N$. We consider 
 ${v}^k\in H^1_0(\O)$  be the unique solution to 
\begin{equation}\label{eq:tvkstf}
 \int_{\O} \n {v}^k\n\f-\int_{\O}\min\left\{b(x),{k}\right\} {v}^k\f=\int_{\O}f\f\quad\forall\f\in C^\infty_c(\O).
\end{equation}
Thanks to   \textbf{Step 1}, we have $v^k\leq u$ in $\O$.\\
Next, we check that  such a sequence  ${v}^k$, satisfying \eqref{eq:tvkstf}, 
  converges to $v$ in $L^2(\O)$ when $k\to\infty$. Indeed, we have 

\begin{eqnarray*}
\| {v}^k\|^2_{{H(\O)} }
&\leq& \| {v}^k\|^2_{{H^1_0(\O)} }-\int_{\O}\min\{b(x),k\}~\!|{v}^k|^2~dx\\
&=&
\int_{\O} f{v}^k~dx\le C\|{v}^k\|_{H(\O)}
\end{eqnarray*}
by H\"older inequality and by (\ref{eq:coerciv}), where the constant $C$
depends on $f$ and $\O$ but not on $k$.
Therefore the sequence ${v}^k$ is bounded in ${H(\O)}$. We conclude that  there exists
$\ti{v}\in {H(\O)}$ such that, for a subsequence, ${v}^k \rightharpoonup \ti{v}$ in
${H(\O)}$.
Now by \eqref{eq:tvkstf},  we have
$$
\la {v}^k,\f\ra_{{H(\O)}}+\int_{\O}\left(b(x)-
\min\{b(x),k\}\right){v}^k\vp=\int_{\O}f\vp.
$$
Since for every $k\geq 1$ and any $\vp\in C^\infty_c(\O) $
 $$\left|\left(b(x)- \min\{b(x),k\}\right){v}^k\vp\right|\leq\left(b(x)- \min\{b(x),k\}\right)u|\vp|
\leq 2b(x)u|\vp|\in L^1(\O),$$
 the dominated convergence theorem implies
that
\begin{equation}\label{eq:vstw}
\la \ti{v},\f\ra_{H(\O)}=\int_{\O} f \vp\quad\textrm{for any $
\vp\in C^\infty_c(\O)$.}
\end{equation}
We therefore have  that $\ti{v}=v$ by uniqueness.
By \eqref{eq:vstw}, we have
\begin{eqnarray*}
 \| {v}-{v}^k\|^2_{{H(\O)}}&=&
\| {v}^k\|_{{H(\O)}}^2- \la {v},{v}^k\ra_{{H(\O)}}
+\la {v},{v}-{v}^k\ra_{{H(\O)}}\\
&=&\| {v}^k\|_{{H(\O)}}^2-\int_{\O} f {v}^k+\la {v},{v}-{v}^k\ra_{{H(\O)}}\\
&\leq&  \| {v}^k\|^2_{H^1_0(\O)}-\int_{\O}\min\{b(x),k\}~\!|{v}^k|^2~dx -\int_{\O} f {v}^k+\la {v},{v}-{v}^k\ra_{{H(\O)}}\\
&=&   \la {v},{v}-{v}^k\ra_{{H(\O)}} .
\end{eqnarray*}
We thus obtain 
$$
C(\O)\,\int_{\O}|{v}-{v}^k|^2~dx\leq \la {v},{v}-{v}^k\ra_{{H(\O)}}\to 0
$$
by \eqref{eq:coerciv}.
Hence  ${v}^k\to {v}$ pointwise and
thus ${v}\leq u$ in $\O$.

 \QED
We conclude this section by pointing out the following Allegretto-Piepenbrink type result which is
essentially contained in \cite{FaMu-ne}. A version for distributional solutions is also contained in [\cite{CFKS}, Theorem 2.12].
\begin{Lemma}\label{lem:AP}
Let $\O$ be a domain (possibly unbounded) in $\R^N$, $N\geq1$. Let $V\in L^1_{loc}(\O)$
and $V>0$ in $\O$. Assume that $u\in L^1_{loc}(\O)$,  $V(x)u\in L^1_{loc}(\O)$  and that $u$ is a
non-negative, non-trivial solution to
$$
-\Delta u\geq V(x) u\quad \mathcal D'(\O).
$$
Then
$$
\int_{\O}|\n \phi|^2~dx\geq \int_{\O}V(x)\,\phi^2~dx \quad\textrm{for any $
\phi\in C^\infty_c(\O)$.}
$$
\end{Lemma}
\proof Put $V_k(x)=\min\{V(x),k\}$ then  Lemma B.1 in  \cite{FaMu-ne} yields
$$
\int_{\O}|\n \phi|^2~dx\geq
\int_{\O}V_k(x)\,\phi^2~dx\quad\textrm{for any $\phi\in C^\infty_c(\O)$.}
$$
To conclude, it suffices to use Fatou's lemma.
 \QED
\begin{Remark} \label{re:ne-cgmu}
Given $\O$ any domain in $\R^N$, $N\geq 1$. Define
$$
\mu(\O):= \inf_{u\in C^\infty_c(\O)}
\frac{\displaystyle\int_{\O}|\n u|^2~dx}
 {\displaystyle\int_{\O}|x|^{-2}u^2~dx}.
 $$
Then  Lemma \ref{lem:AP} clearly implies that
  if $ c> \mu(\O)$
there is no non-negative and non-trivial $u\in L^1_{loc}(\O)$ that
 satisfies $ -\Delta u-\frac{c}{|x|^2} u\geq0$ in $\calD'(\O)$.
\end{Remark}
%
Suppose that $\O$ is a smooth bounded domain and  that  the potential $b(x)$ satisfies
$$
 \int_{\O}|\n \vp|^2~dx-\int_{\O}b(x)\vp^2~dx\geq C(b)\left( \int_{\O}|\vp|^r~dx\right)^{\frac{2}{r}}\quad\textrm{for any
 $\vp\in C^\infty_c(\O)$}
$$
for some $C(b)>0$ and  $2<r$. By [\cite{DD} Lemma 7.2], we can let $G\in L^1(\O\times\O)$ be the Green function associated to $-\D -b(x)$:
$$
\begin{cases}
-\D G(\cdot,y)-b(x) G(\cdot ,y)=\d_y \quad\textrm{ in } \O,\\
G(\cdot,y)=0\quad\textrm{ on } \de\O,
\end{cases}
$$
where $\d_y$ denotes the Dirac measure at some $y\in\O$. Define
$$
\zeta_0(x):=\int_\O G(x,y)\, dy
$$
which is the $H(\O)$-solution  to $-\D\zeta_0 -b(x)\,\zeta_0=1$.
By using Lemma \ref{lem:min-sol} and Lemma \ref{lem:AP}, we can prove the following
\begin{Proposition}
Suppose that  $\displaystyle\int_{\O} \zeta_0^{p+1}\,dx =\infty$ for some $p>r$ then there is
no nonnegative and nontrivial $u$ satisfying $-\D u -b(x)\,u\geq u^p$ in $\calD'(\O)$.
\end{Proposition}
\proof
If such $u$ exists, it is positive by the maximum principle  therefore, we can define $v\in H(\O)$ be the solution
of $-\D v -b(x)\,v=\min(u^p,1)$ so that by Lemma \ref{lem:min-sol} we have $u\geq v$ in $\O$.
 Thanks to [\cite{DD} Corollary 2.4],
we have $u\geq v\geq C \zeta_0$. By applying  Lemma \ref{lem:AP} with $V(x)=b(x)+(C\zeta_0)^{p-1}$ we conclude that
$$
\infty>\|\zeta_0\|_{H(\O)}\geq C^{p+1} \int_{\O} \zeta_0^{p+1}\,dx.
$$
\QED

\section{Proof of Theorem \ref{th:neC}}
\label{S:cone}
We state the following lemma which is a consequence  of Lemma \ref{lem:min-sol}
and [\cite{LLM}, Theorem 4.2].
\begin{Lemma}\label{lem:est-cone}
Let $u\in L^1_{loc}(\O)$ be  positive  and let $f\in L^1_{loc}(\calC_{\Sig}^r)$ with $f\gneqq0$ such that 
$$
-\D u-V\left(\frac{x}{|x|}\right)\,|x|^{-2}\,u\geq f \quad\textit{ in $ \calD'(\calC_{{\Sig}}^{r})$,}
$$
where $\|V\|_{ L^\infty(\Sig)}\leq \mu(\calC_{\Sig})$ and $V\geq0$.
Then for every $ \ti{\Sig}\subset\subset {\Sig}$  there exists a  constant $C>0$ such that
\be\label{eq:estofu}
u(x)\geq C |x|^{\frac{2-N}{2}+\sqrt{ \frac{(2-N)^2}{4}+\l_{1,V}   }} \quad
 \textrm{ in }\calC_{\ti{\Sig}}^{r/2},
\ee
where   $\l_{1,V}$  is the first Dirichlet eigenvalue of  
 $-\D_{\S^{N-1}}\Phi- V \Phi= \l_{1,V}\Phi$ on $ \Sig$ . 
\end{Lemma}
\proof
Up to a scaling, we can assume that $r=1$. We recall the following
improved Hardy inequality
\begin{equation}
\label{eq:Leray} \int_{\C^1_\Sigma}|\nabla \vp|^2~dx-\mu(\C_\Sigma)
\int_{\C^1_\Sigma}|x|^{-2}|\vp|^2\\
\ge
C_0\int_{\C^1_\Sigma}|\vp|^2~dx~\!\quad\forall \vp\in C^\infty_c(\C^1_\Sig),
\end{equation}
for some $C_0>0$  (see for instance \cite{FaMu1}).  We can therefore pick 
$v\in H\left(\C_\Sigma^1\right)$ solves
\begin{equation}
\label{eq:v}
-\Delta v- V\left(\frac{x}{|x|}\right)\, |x|^{-2}v=~\!\min(f,1) \quad\textrm{in $H\left(\C_\Sigma^1\right)$.}
\end{equation}
Then  by the maximum principle and  Lemma \ref{lem:min-sol},
 we have $0<v\le u$ in ${\C_\Sigma^1}$.\\
Approximating $v$ by smooth functions  compactly supported in $ \C_\Sigma^1$  with respect to the 
 $ H\left(\C_\Sigma^1\right)$-norm, we infer that
$$
-\Delta v- V\left(\frac{x}{|x|}\right)\, |x|^{-2}v=~\!\min(f,1) \quad\textrm{in $\calD'\left(\C_\Sigma^1\right)$.}
$$
Elliptic regularity theory then 
implies that $v\in C^{1,\g}_{loc}(\C_\Sigma^1)\subset H^{1}_{loc}(\C_\Sigma^1)$.
By applying [\cite{LLM}, Theorem 4.2] (up to Kelvin transform), we get the lower estimate 
\eqref{eq:estofu} for $v$ and hence for $u$.

\QED

\bigskip\noindent
\textbf{Proof of Theorem \ref{th:neC}}\\
\noindent
Up to a scaling, we can  assume that $R=1$. We argue by contradiction. If  $u\neq0$ then by the maximum principle
  $u>0$ in ${\C_\Sigma^1}$.
We will show that   appropriate lower bound of $u$ 
and an application of Lemma \ref{lem:AP} will lead to a contradiction.\\

\bigskip
\noindent
{\bf Case 1:  $c<\mu(\calC_\Sig)$.}\\
By Lemma \ref{lem:est-cone}
$$
u(x)\geq C_0|x|^{\frac{2-N}{2}+\sqrt{\mu(\calC_\Sig)-c}}\quad \forall x\in  \calC_{\ti{\Sig}}^{1/2},
$$
where $C_0$ is a positive constant and $\ti{\Sig}\subset\subset\Sig $. 
By assumption $u^{p-1}(x)|x|^2\geq C_0^{p-1}$. 
In particular  for every $\e\in(0,1)$, we have
$$
-\D u-(c+\e V)|x|^{-2}u\geq \frac{1}{2}u^{p} \quad \textrm{ in } \calD'( \C_\Sigma^{1/2}),
$$
where $V= \frac{C_0^{p-1}}{2} \chi_{  {\ti{\Sig}}} $. 
We notice that for $\e$ small, $c+\e V< \mu(\calC_\Sig)$. 
 We apply once more  Lemma \ref{lem:est-cone} to get
\be\label{eq:estveps}
u(x)\geq  C_1|x|^{\frac{2-N}{2}+\sqrt{(N-2)^2/4+\l_{1,\e}}}\quad \forall x\in  \calC_{\ti{\Sig}}^{1/4},
\ee
where $ \l_{1,\e}$ is the first Dirichlet  eigenvalue of
 $-\D_{\S^{N-1}}\Phi-(c+\e V )\Phi= \l_{1,\e}\Phi$ on $ \Sig$.
We observe that, for $\e$ small, $\l_{1,\e}<\l_1(\Sig)-c<0$ and thus
$$
p-1\geq \frac{2}{\a_{\Sig}^-}>\frac{2}{ \frac{N-2}{2}-\sqrt{(N-2)^2/4+\l_{1,\e}}   }>0.
$$
Recalling that $-\D u\geq u^{p-1}u$, we deduce  from \eqref{eq:estveps} that 
$$
-\D u-\rho(x)  |x|^{-2} u\geq0  \quad \textrm{ in } \calD'( \C_{\ti\Sigma}^{1/4}),
$$
where $ \rho(x)=\frac{C_1^{p-1}}{2}|x|^{(\frac{2-N}{2}+\sqrt{(N-2)^2/4+\l_{1,\e}})(p-1)+2}$.
 Since $\rho(x)\to+\infty$ as $|x|\to 0$,
 applying  Lemma \ref{lem:AP}, we contradict the sharpness
 of the Hardy constant $\mu(\calC_{\ti{\Sig}})$.

\bigskip
\noindent
{\bf Case 2:  $c=\mu(\calC_\Sig)$.}\\
We consider the function  $v\in H(\C_{\Sig}^1)$ solving 
 $$-\D v-\mu(\C_{\Sig} )|x|^{-2} v=\min(u^p,1).$$
Then by Lemma \ref{lem:min-sol} and the maximum principle $0<v\leq u$ in $ \C_{\Sig}^1$.
By Lemma \ref{lem:est-cone}, 
$$
v(x)\geq C |x|^{\frac{2-N}{2}}\qquad \textrm{ for } x\in \calC_{\ti{\Sig}}^{1/2}.
$$
Since $-\D u - \mu(\C_{\Sig} )|x|^{-2} u=u^{p-1}u$ in $\calD'( \C_{\Sig}^1)$, by Lemma \ref{lem:AP} and the above estimate, we have
\begin{eqnarray*}
\infty>\|v\|_{H(\C_{\Sig}^1)}^2&\geq&
 \int_{\C_{\Sig}^{1}}
u^{p-1}v^2\,dx\geq \int_{\C_{\Sig}^{1}}
v^{p+1}\,dx\\
&\geq& C\int_{\C_{\ti{\Sig}}^{1/2}}
v^{p+1}\,dx\geq C\int_{\C_{\ti{\Sig}}^{1/2}}
|x|^{-N }\,dx=\infty.
\end{eqnarray*}
This readily leads to a contradiction.
Theorem \ref{th:neC} is completely proved.
 \QED
 %
\section{Smooth  domains}
\label{S:smooth}
In this section, we introduce  a system of   coordinates near $0\in\de\O$ that flattens  $\de\O$, see \cite{MMF-PJM}.
This will allows us to  construct a (super-)
sub-solution via the function $y^1|y|^{-\frac{N}{2}+\sqrt{\frac{N^2}{4}}-c}$ which is
 the (virtual) ground state for the operator $\D +c|y|^{-2}$ in the half-space $\R^{N}_+$.
\subsection{Fermi coordinates}\label{s:pn} 
We denote by  $\{E_1,E_2,\dots,E_N\}$ the standard orthonormal basis
of $\R^N$ and we put
\begin{gather*}
\R^N_+=\{y\in\R^N\,:\,y^1>0\}~\!,\quad \S^{N-1}_+=\S^{N-1}\cap \R^N_+~\!,\\
B_r(y_0)=\{y\in\R^N\,:\,|y-y_0|<r\}~,\quad
B_r^+=B_r(0)\cap\R^N_+~\!.
\end{gather*}
Let $\calU$ be an open subset of $\R^N$ with boundary $\M:=\de\calU$
a smooth closed hyper-surface of ${\R^N}$ and $0\in\M$. We write
$N_{\M} $ for the unit normal vector-field of $\M$ pointed into
$\calU$. Up to a rotation, we assume that $N_\M(0)=E_1$.  For
$x\in\R^N$, we let $d_\M(x)=\textrm{dist}(\M,x)$ be the distance
function of $\M$. Given $x\in\calU$  and close to $\M$ then it can
be written uniquely as $x=\s_x+d_{\M}(x)\,N_{\M}(\s_x)$, where
$\s_x$ is the projection of $x$ on $\M$. We further   use the
\textrm{Fermi coordinates} $(y^2,\dots,y^N)$ on $\M$ so that for
$\s_x$ close to 0, we have
$$
\s_x=\textrm{Exp}_0\left(\sum_{i=2}^N y^iE_i\right),
$$
where $\textrm{Exp}_0:\R^{N-1}\to \M$ is the exponential mapping on $\M$
endowed with the metric induced by $\R^N$, see  \cite{Docarmo}. In this way a
neighborhood of 0 in $\calU$ can be parameterize by the map
$$
 F_{\M}(y)=\textrm{Exp}_0\left(\sum_{i=2}^Ny^iE_i\right)+y^1\,N_{\M}\left(
\textrm{Exp}_0\left(\sum_{i=2}^Ny^iE_i\right)\right),\quad y\in
B_r^+,
$$
for some $r>0$. In this coordinates, the Laplacian $\D$  is given
by
$$
\sum_{i=1}^N \frac{\de^2 }{(\de x^i)^2} =\frac{\de^2 }{(\de y^1)^2}+h_{\M}\circ F_\M\,\frac{\de }{\de
y^1}+\sum_{i,j=2}^N\frac{\de}{\de y^i}\left(\sqrt{|g|}{g^{ij}}
\frac{\de }{\de y^j} \right),
$$
where $h_{\M}(x)=\D\, d_{\M}(x)$; for $i,j=2\dots,N$,
$g_{ij}=\langle\frac{\de F_\M }{\de y^i},\frac{\de F_\M}{\de y^j}
\rangle$; the quantity ${|g|}$ is the determinant of $g$ and
$g^{ij}$ is the component of the
inverse of the matrix $(g_{ij})_{2\leq i,j\leq N}$.\\
Since $g_{ij}=\d_{ij}+O(y^1)+O(|y|^2)$ (see \cite{MMF-PJM}), we have the following
Taylor expansion
\be\label{eq:expDFc}
\sum_{i=1}^N \frac{\de^2 }{(\de x^i)^2} =\sum_{i=1}^N \frac{\de^2 }{(\de y^i)^2}+h_{\M}\circ F_\M\,\frac{\de }{\de
y^1}+\sum_{i=2}^NO_i(|y|)\frac{\de }{\de y^i}+\sum_{i,j=2}^N
O_{ij}(|y|)\frac{\de^2 }{\de y^i\de y^j}.
\ee
For $a\in\R$, we put  $X_a(t):=|\log t|^a$, $t\in(0,1)$ and for $c\leq{\frac{N^2}{4}}$, set
$$
\ov{\o}_a(y):=y^1|y|^{-\frac{N}{2}+\sqrt{{\frac{N^2}{4}}-c}}X_a(|y|)\quad
\forall y\in\R^N_+
$$
and put
$$
L_y:=- \sum_{i=1}^N \frac{\de^2 }{(\de y^i)^2} -c\,|y|^{-2}+a(a-1)|y|^{-2}X_{-2}(|y|).
$$
 Then  one easily
verifies that
$$
\begin{cases}
\displaystyle
L_y\,\ov{\o}_a=2a\sqrt{{\frac{N^2}{4}}-c}\,|y|^{-2}X_{-1}(|y|)\,\ov{\o}_a&
\textrm{ in $\R^N_+$,}\\
\displaystyle\ov{\o}_a=0&\textrm{ on $\de \R^N_+\setminus\{0\}$,}\\
\displaystyle\ov{\o}_a\in H^1(B_R^+)\,\, \forall R>0,\,\,\forall
c<{\frac{N^2}{4}},\,\,\forall a\leq 0.
\end{cases}
$$
For $K\in\R$, we define $$\o_{a,K}(y)=e^{K y^1}\,\ov{\o}_a(y).$$
This function satisfies similar boundary and integrability
conditions as $\ov{\o}_a$. In addition  it holds that
%
\be\label{eq:Lyoak}
\begin{array}{c}
\displaystyle L_y\,\o_{a,K} =-\frac{2K}{y^1}\o_{a,K}+ \displaystyle
2a\sqrt{{\frac{N^2}{4}}-c}\,|y|^{-2}X_{-1}(|y|)\,\ov{\o}_a \displaystyle\\
\hspace{5cm}+2K
\left(\frac{N}{2}-\sqrt{{\frac{N^2}{4}}-c}+aX_{-1}(|y|)\right)\frac{y^1}{|y|^2}\o_{a,K}-K^2\o_{a,K}.
\end{array}
\ee
Furthermore   for all $a\in\R$
\begin{eqnarray*}
\sum_{i=2}^NO_i(|y|)\frac{\de \o_{a,K}}{\de y^i}+
\sum_{i,j=2}^NO_{ij}(|y|)\frac{\de^2 \o_{a,K}}{\de y^i\de
y^j}&=&y^1\,e^{Ky^1}
O\left(|y|^{-\frac{N}{2}-1+\sqrt{{\frac{N^2}{4}}-c}}X_a(|y|)\right)\\
&=&\calO_{a,K}(|y|^{-1})\,\o_{a,K}(y).
\end{eqnarray*}
Here the error term $\calO_{a,K}$  has the property that for any
$A>0$ and $c_0<{\frac{N^2}{4}}$, there exit some  constants $C>0$ and $s_0>0$  such that
\be\label{eq:errterm}
\left|\calO_{a,K}(s)\right|\leq C\,s\quad\forall
s\in(0,s_0),\,\,\forall a\in[-A,A]\,\,\forall c\in\left[c_0,{\frac{N^2}{4}}\right].
\ee
Let
\be\label{eq:WaK}
W_{a,K}(x):=\o_{a,K}(F_{\M}^{-1}(x)),\qquad\forall\,
x\in\calB_r^+:=F_\M(B_r^+).
\ee
 Then using \eqref{eq:expDFc}, \eqref{eq:Lyoak} and the
fact that $|x|=|y|+O(|y|^2)$ we obtain the following expansions
\be\label{eq:expD}
L_x\,W_{a,K}=-\left(\frac{2K+h_{\M}(x)}{d_{\M}(x)}\right)\,W_{a,K}+
2a\sqrt{{\frac{N^2}{4}}-c}\,|x|^{-2}X_{-1}(|x|)\,W_{a,K}
+\calO_{a,K}(|x|^{-1})\,W_{a,K},
\ee
with $ L_x:=-\D -c\,|x|^{-2}+a(a-1)|x|^{-2}X_{-2}(|x|).$
 Moreover it is easy to see that
\be
 \begin{cases}
\label{eq:wH1}
W_{a,K}>0&\textrm{ in $\calB_r^+$,}\\
W_{a,K}=0&\textrm{ on $\M\cap
\de\calB_r^+\setminus\{0\}$,}\\
W_{a,K}\in H^1(\calB_r^+)\quad\,\forall c<{\frac{N^2}{4}},\, \textrm{ and }\,\,\forall a\leq 0.
\end{cases}
\ee
\subsection{Non-existence}
We start by recalling the following local improved Hardy inequality.
 Given a  domain $\O\subset\R^N$, of class $C^2$
at $0\in\de\O$, there exist two
 constants $C(\O)>0$ and $r_0=r_0(\O)>0$ such that
\begin{equation}\label{eq:impH}
\int_{\O_{r_0}}|\n u|^2\,dx-c\int_{\O_{r_0}}|x|^{-2}u^2\,dx\geq
C(\O)\,\int_{\O_{r_0}}u^2\,dx\quad\forall u\in C^\infty_c(\O_{r_0}),
\end{equation}
for every $c\in\left(-\infty,{\frac{N^2}{4}}\right]$, with $\O_{r_0}:=\O\cap B_{r_0}(0)$,
 see \cite{mmf}. From this we can define the
space $H(\O_{r_0})$  to be the completion of $C^\infty_c(\O_{r_0})$
with respect to the
 scalar product
 $$
\int_{\O_{r_0}}\n u\n v-c\int_{\O_{r_0}}|x|^{-2}uv\quad\forall
u,v\in C^\infty_c(\O_{r_0}).
 $$
 In the
 sequel we  will  assume that $\O$ contains the ball $B=B_1(E_1)$ such
that $\de{B}\cap\de{\O}=\{0\}$.
 Recalling the  notations in Section \ref{s:pn}, we state the following result
\begin{Lemma}\label{lem:estbl}
 Let $c_0\in\left(-\infty,{\frac{N^2}{4}}\right]$ and $f\in L^\infty(\O_{r_0})$ be a
non-negative and non-trivial function. For $c\in \left[c_0,{\frac{N^2}{4}}\right]$, let $v\in H(\O_{r_0})$ be
the unique solution of the problem
$$
\int_{\O_{r_0}}\n
v\n\phi\,dx-c\int_{\O_{r_0}}|x|^{-2}v\phi\,dx=\int_{\O_{r_0}}f\phi\,dx\quad
\forall\phi\in C^\infty_c(\O_{r_0}).
$$
Then    there exist $R>0$ and $r>0$  such that
$$
v(F_{\de B}(y))\geq
R\,y^1|y|^{-\frac{N}{2}+\sqrt{{\frac{N^2}{4}}-c}} \quad\forall y\in
B_{r}^+\quad\forall c\in\left[c_0,{\frac{N^2}{4}}\right].
$$
\end{Lemma}
\proof
 For  $a\leq 0$ and $r>0$  small we define $G_r^+:=F_{\de B}(B_r^+)$ and
$$
w_{a}(x):=\o_{a,N-1}(F_{\de B}^{-1}(x)),\quad \forall x\in G_r^+.
$$
 Letting
$L:=-\D-c\,|x|^{-2}$, by
 \eqref{eq:expD},
\begin{eqnarray*}
L\,(w_0+w_{-1})&\leq&  -  \frac{ 2(N-1)+ h_{\de B} }{d_{{\de B}}}\,{(w_0+w_{-1})}
-2\,|x|^{-2}\,X_{-2}(|x|)
\,w_{-1}\\
&-&2\sqrt{{\frac{N^2}{4}}-c}|x|^{-2}X_{-1}(|x|)\,w_{-1}+O(|x|^{-1})\,(w_0+w_{-1}).
\end{eqnarray*}
We observe that
\be\label{eq:wm1pw0}
 w_{-1}(x)=w_0(x) |\log|F_{\de B}^{-1}(x) ||^{-1} = w_0(x)\left( X_{-1}(|x|)+ O(|x|)\right)  .
\ee
Since   $-h_{\de B}(x)=(N-1)\left(1+O(|x|)\right)$ in $G_r^+$, we have using 
   \eqref{eq:wm1pw0},
\begin{equation}\label{eq:lwaneg}
L\,(w_0+w_{-1})\leq0\quad\textrm{ in } G_r^+,
\end{equation}
for  $r$ positive small.\\
\textbf{Case $c\in\left[c_0,{\frac{N^2}{4}}\right)$}.\\
We put $U=w_0+w_{-1}$. Then  $U\in H^1(G_r^+)\cap C(G_r^+)$ by \eqref{eq:wH1} and $v\in
H^1_0(\O_{r_0})\cap C(\O_{r_0})$
 by elliptic regularity theory and Remark \ref{R:H2}. Moreover $v>0$ in  $\O_{r_0}$ by 
the maximum principle.
Therefore since ${F_{\de B}\left(\ov{r \S^{N-1}_+}\right)} \subset
\O_{r_0}$, we can let
\begin{equation}\label{eq:choiR}
R=r^{\frac{N-2}{2}}\, e^{-(N-1)r}\,\inf_{y\in \ov{r\S^{N-1}_+}}\,v>0
\end{equation}
 so that
$$
R\,{U}\leq v\quad\textrm{ on } {F_{\de B}\left(\ov{r
\S^{N-1}_+}\right)}.
$$
By \eqref{eq:wH1} and setting $\vp=R\, {U}-v$, we get 
$\vp^+:=\max(\vp,0)\in H^1_0(G_r^+)$ because $U=0$ on $\de
B\cap\de G_r^+$. Since $Lv\geq0$, we have
$$
L\,\vp\leq 0\quad\textrm{ in }G_r^+,
$$
by  \eqref{eq:lwaneg}.
 Multiplying the above inequality by $\vp^+$ and integrating by parts
 yields
$$
\int_{G_r^+}|\n \vp^+|^2\,dx-c\int_{G_r^+}|x|^{-2}|\vp^+|^2\,dx
\leq0.
$$
This  implies that $\vp^+\equiv0$ in
$G_r^+$ for all $r$ positive small. We conclude that $v\geq R( w_0+w_{-1})$ in $ G_r^+$ and thus
$$
v(F_{\de B}(y))\geq R\,w_0(F_{\de B}(y))\geq R\, y^1|y|^{-\frac
{N}{2}+\sqrt{{\frac{N^2}{4}}-c}} ,\quad\forall y\in G_{r}^+,\quad\forall
c\in\left[c_0,{\frac{N^2}{4}}\right).
$$
\textbf{Case $c={\frac{N^2}{4}}$}.\\
In this case, we notice that the solutions $v_k$ to the problem
$$
\int_{\O_{r_0}}\n
v_k\n\phi\,dx-\left({\frac{N^2}{4}}-\frac{1}{k}\right)\int_{\O_{r_0}}|x|^{-2}v_k\phi\,dx
=\int_{\O_{r_0}}f\phi\,dx\quad \forall\phi\in H(\O_{r_0})
$$
are $H^1_0$-solutions if $r_0$ is small enough (independent on $k$) and
 they are monotone increasing to $v$ as
$k\to \infty$.  Hence by \eqref{eq:choiR} and  from the above argument we
deduce that there exist  an
integer $k_0\geq1$ and a constant $R$ (possibly depending on $k_0$)
such that
$$
v(F_{\de B}(y))\geq v_{k}(F_{\de B}(y))\geq R\, y^1|y|^{-\frac
{N}{2}+\sqrt{\frac{1}{k}}} ,\quad\forall y\in G_{r}^+,\quad\forall k\geq k_0.
$$
Passing to the limit as $k\to \infty$, we get the result.
 \QED
%
%
\bigskip
\noindent
\subsubsection{{Proof of Theorem \ref{th:ne-sm-smoo-loc}}}
Recall that
\be\label{eq:assumpp}
\alpha^-_{\S^{N-1}_+}= \frac{N-2}{2}-\sqrt{\frac{N^2}{4}-c },\qquad p\geq p_{\S^{N-1}_+}:=1+\frac{2}{\alpha^-_{\S^{N-1}_+}},
\ee 
Suppose that $u\neq0$ near 0 thus we can  find a bounded function 
$f$ with $f\neq0$ and $0\leq f\leq u^p$. By Lemma
\ref{lem:min-sol} and the maximum principle, there exits $v\in H(\O_{r_0})$ such that $u\geq v>0$
and
$$
\int_{\O_{r_0}}\n
v\n\phi\,dx-c\int_{\O_{r_0}}|x|^{-2}v\phi\,dx=\int_{\O_{r_0}}f\phi\,dx\quad
\forall\phi\in C^\infty_c(\O_{r_0}),
$$
for some $r_0>0$ small. In addition 
 Lemma \ref{lem:estbl}  yields
\begin{equation}\label{eq:vfgmy1}
v(F_{\de B}(y))\geq R\,{y^1}
|y|^{-\frac{N}{2}+\sqrt{{\frac{N^2}{4}}-c}}  \quad\forall
y\in\,B_r^+.
\end{equation}

\bigskip
\noindent
\textbf{Case 1: $c\in\left(N-1,{\frac{N^2}{4}}\right)$}.\\
Since $ -\frac{N}{2}+\sqrt{{\frac{N^2}{4}}-c}<0$, \eqref{eq:vfgmy1} implies that 
\begin{equation}\label{eq:vfgmynp}
u(x)\geq v(x)\geq C\,{d_{\de B}(x)}\,|x|^{-\frac{N}{2}+\sqrt{{\frac{N^2}{4}}-c}} 
      \quad\forall
x\in\,G_r^+,
\end{equation}
where we recall that $d_{\de B}(F_{\de B}(y))=y^1$  and $|F_{\de B}^{-1}(x)|\leq C|x|$. \\
Let $\g\in(0,1)$
then for every $x\in B_{\g}(\g E_1)\subset B_1(E_1)=B$, we have
$$
d_{\de B}(x)=1-|x-E_1|> (1-\g)\,x^1.
$$
Using this together with  \eqref{eq:assumpp} and \eqref{eq:vfgmynp}, we obtain
\begin{equation}\label{eq:vfgmy}
u^{p-1}(x)\geq C\,\left(\frac{x^1}{|x|}\right)^{p-1}\,|x|^{-2 } 
     \quad\forall
x\in\,G_r^+.
\end{equation}
Since for $\g>0$ small $u$ satisfies
$$
-\D u-c|x|^{-2}u\geq \frac{1}{2}u^{p-1} u+\frac{1}{2}u^{p} \quad \textrm{ in } \calD'(B_{\g}(\g E_1)),
$$
we thus have from \eqref{eq:vfgmy}
$$
-\D u-V_\e\left(\frac{x}{|x|}\right)|x|^{-2}u\geq \frac{1}{2}u^{p} \quad \textrm{ in } \calD'(B_{\g}(\g E_1)),
$$
where $V_\e\left(\frac{x}{|x|}\right)=c+ \frac{\e}{2} C\,\left(\frac{x}{|x|}\cdot E_1\right)^{p-1}$ 
for every $\e\in (0,1)$. From now on, we will fix $\e$ so small that $V_\e<\frac{N^2}{4}$.\\
Given $\d\in(0,1)$, consider the cone $\calC_\d:=\{x\in\R^N\,:\,x^1>\d|x|\}$
 and define $\Sig_\d=\calC_\d\cap\S^{N-1}$. We observe that for every $\d\in(0,1)$, 
there exists $r_\d>0$ such that the cone-like domain
 $$\calC_{\Sig_\d}^{r_\d}\subset B_{\g}(\g E_1) .$$
It follows that  
$$
-\D u-V_\e\left(\frac{x}{|x|}\right)|x|^{-2}u \geq \frac{u^p}{2} \quad \textrm{ in } \calD'(\calC_{\Sig_\d}^{r_\d}).
$$
Let $ \l_{1,\d, \e}$ be the
 first Dirichlet  eigenvalue of $-\D_{\S^{N-1}}\Phi- V_\e \Phi= \l_{1,\d,\e}\Phi$ on
 $\Sig_\d$.
Since $\l_1(\Sig_\d)\searrow N-1=\l_1(\S^{N-1}_+)$ as $\d\to0$, 
we can choose a $\d_\e \in(0,1)$ 
such that    
\be
\label{eq:upl1de}
\l_{1,\d,\e}<N-1-c<0\quad \forall \d\in(0,\d_\e).
\ee
Since $V_\e\leq \frac{N^2}{4}<\frac{(N-2)^2}{4}+\l_1(\Sig_\d)=\mu(\calC_{\Sig_\d})$
 for every $\d\in(0,\d_\e)$, we can apply  Lemma \ref{lem:est-cone} to
 have $\forall \d\in(0,\d_\e)$
$$
 u(x)\geq C |x|^{\frac{2-N}{2}+\sqrt{ \frac{(2-N)^2}{4}   +\l_{1,\d,\e}}},\quad \textrm{ in }
\calC_{\ti{\Sig}_\d}^{r_\d/2}
$$
where $\ti{\Sig}_\d\subset\subset\Sig_\d$.
We  get from \eqref{eq:upl1de}
$$
p-1\geq \frac{2}{\a_{\S^{N-1}_+}^-}>\frac{2}{ \frac{N-2}{2}-\sqrt{(N-2)^2/4+\l_{1,\d, \e}}   }>0.
$$
 Since $-\D u\geq u^{p-1}u$, we deduce   that $\forall \d\in(0,\d_\e)$
$$
-\D u-\rho(x)  |x|^{-2} u\geq0  \quad \textrm{ in } \calD'( \calC_{\ti{\Sig}_\d}^{r_\d/2} ),
$$
where $ \rho(x)\geq C |x|^{(\frac{2-N}{2}+\sqrt{(N-2)^2/4+\l_{1,\d,\e}})(p-1)+2}$. 
Since $\rho(x)\to+\infty$ as $|x|\to 0$,
 applying  Lemma \ref{lem:AP}, we contradict the sharpness of the Hardy 
constant $\mu(\calC_{\ti{\Sig}_\d})$.

\bigskip
\noindent
\textbf{Case 2: $c={\frac{N^2}{4}}$}.\\
 Here, we recall that   $p\geq \frac{N+2}{N-2}$. By  \eqref{eq:impH} we can let $\zeta\in H(G_r^+)$ be the unique
 solution to the problem
 $$
\int_{G_r^+}\n
\zeta\n\phi\,dx-{\frac{N^2}{4}}\int_{G_r^+}|x|^{-2}\zeta\phi\,dx
=\int_{G_r^+}1\,\phi\,dx\quad
\forall\phi\in C^\infty_c(G_r^+).
$$
We put  $\Phi(y)=\frac{y^1}{|y|}$. Then by  Lemma \ref{lem:AP},
 Lemma \ref{lem:estbl} and \eqref{eq:vfgmy1}, we have
\begin{eqnarray*}
\|\zeta\|_{H(G_r^+)}^2&\geq& \int_{G_r^+}v^{p-1}|\zeta|^2\,dx\\
&\geq& C\int_{B_r^+}
|y|^{\frac{2-N}{2}(p+1)} \Phi^{p+1}\,\sqrt{|\hat{g}|}(y)dy\\
&\geq& C \int_{B_r^+}
|y|^{\frac{2-N}{2}(p+1)} \Phi^{p+1}\,dy\\
&=& C \int_{\S^{N-1}_+}\Phi^{p+1}\,d\s\int_0^r
t^{\frac{2-N}{2}(p+1)} t^{N-1}\,dt\\
 &\geq& C
\int_{\S^{N-1}_+}\Phi^{p+1}\,d\s\int_0^r t^{-1} \,dt=\infty.
\end{eqnarray*}
This clearly contradicts the fact that $\zeta\in H(G_r^+)$.
 \QED

\subsection{Existence}\label{s:exits}
Let $\O$ be a  domain of $\R^N$, $N\geq3$ which is of class $C^2$ at $0\in\O$, we shall show that for some $r>0$ small, there exists a positive function
$u\in L^p(\O\cap B_r(0))$, $ 1<p< p_{\S^{N-1}_+}=1+\frac{2}{\a^-_{\S^{N-1}_+}}$ and
$$
-\D u-\frac{c}{|x|^{2}}u\geq u^p\quad\textrm{ in } \calD'(\O\cap B_r(0)).
$$
Letting $B$ be a unit ball with $0\in\de B$, call $\calU=\R^N\setminus \ov{B}$ and $\calM=\de\calU$. Under the notations in Section \ref{s:pn},
the above existence result is  a consequence of the following
\begin{Proposition}\label{prop:locex}
Let $1<p<p_{\S^{N-1}_+}$ and $N-1< c \leq \frac{N^2}{4}$. Then there exists $r>0$ small such that the problem
\be
\begin{cases}
-\D w-\frac{c}{|x|^{2}}w= w^p\quad\textrm{ in } \calD'(\calB_r^+), \\
w\in L^p(\calB_r^+ ),\\
w>0\quad \calB_r^+
\end{cases}
\ee
has a    supersolution, with $\calB_r^+=F_{\calM}(B_r^+)$.
\end{Proposition}
\proof
Notice that  $h_{\M}(x)=\frac{N-1}{1+d_\M(x)}$ and
thus
\be\label{eq:ekphgnm1}
-\frac{2(1-N)+h _{\M}(x)}{d_\M(x)}\geq
\frac{N-1}{d_\M(x)}\quad\forall x\in \calU.
\ee
Define (see \eqref{eq:WaK})
$$
w(x)=\o_{\frac{1}{2p},1-N}(F^{-1}_{\calM}(x))\quad\forall x\in \calB_r^+.
$$
By \eqref{eq:expD}, \eqref{eq:ekphgnm1} and using the fact that $|x|= |y|+O(|y|^2)$, we have
\begin{eqnarray*}
-\D w (F_{\calM}(y)) -c |F_{\calM}(y)|^{-2}w(F_{\calM}(y))&\geq&\frac{2p-1}{4 p^2}|y|^{-2}X_{-2}(|y|)\,w( F_{\calM}(y))\\
&&\,+O(|y|^{-1})\,w(F_{\calM}(y)).
\end{eqnarray*}
In particular if $r>0$ small
\begin{eqnarray*}
-\D w (F_{\calM}(y)) -c |F_{\calM}(y)|^{-2}w(F_{\calM}(y))&\geq& C|y|^{-2}X_{-2}(|y|)\,w( F_{\calM}(y))\quad \forall y\in B_r^+,
\end{eqnarray*}
with $C>0$ a constant depending only on $p$ and $N$.
Therefore $w$ is a supersolution provided
$$
C|y|^{-2}X_{-2}(|y|)\,w( F_{\calM}(y))\geq w(F_{\calM}(y))^p\quad \forall y\in B_r^+
$$
or equivalently
$$
C \frac{y^1}{|y|}e^{(1-N)y^1}|y|^{-\a^-_{\S^{N-1}_+}-2}\left|\log|y|\right|^{-2+\frac{1}{2p}}
\geq
\left(\frac{y^1}{|y|}e^{(1-N)y^1}|y|^{-\a^-_{\S^{N-1}_+}}\left|\log|y|\right|^{\frac{1}{2p}} \right)^p\quad \forall y\in B_r^+.
$$
Since $0<\frac{y^1}{|y|}e^{(1-N)y^1}\leq 1 $ and $p\geq 1$, the above holds if
$$
C|y|^{-\a^-_{\S^{N-1}_+}-2}\left|\log|y|\right|^{-2+\frac{1}{2p}}
\geq
\left(|y|^{-\a^-_{\S^{N-1}_+}}\left|\log|y|\right|^{\frac{1}{2p}} \right)^p\quad \forall y\in B_r^+.
$$
The previous inequality is true provided
$$
C|y|^{-\a^-_{\S^{N-1}_+}-2+p\a^-_{\S^{N-1}_+}} \left|\log|y|\right|^{-2+\frac{1}{2p}-\frac{1}{2}}\geq1\quad \forall y\in B_r^+.
$$
This is clearly possible  whenever $p< p_{\S^{N-1}_+}=1+\frac{2}{\a^-_{\S^{N-1}_+}}$ and $r>0$ is small enough. Finally, we notice that
$
\displaystyle\int_{\calB_r^+}w^p\,dx\leq C \displaystyle\int_0^rt^{\frac{N-4}{2}-\sqrt{\frac{N^2}{4}-c}}|\log t|^{\frac{1}{2}}\,dt<\infty
$, when $N-1<c\leq \frac{N^2}{4}$. This concludes the proof.
\QED

\section{Problem with perturbation}\label{s:pblmp}
We let ${\Gamma}\subset\R^N$ be a {smooth closed submanifold} of {dimension $k$ with $1\leq k< N-2$}.
Let $\O$ be a smooth domain in $\R^{N}$ containing $\G$.
We study the problem
\be\label{eq:petplmi}
\begin{cases}
-\D u - \frac{(N-k-2)^2}{4}\,\frac{1}{\textrm{dist}(x, \Gamma)^{2}}\, {q(x)}\,u\geq u^p\quad\textrm{ in } \calD'(\O\setminus \Gamma), \\
u\in L^p_{loc}(\O\setminus \Gamma),\\
u\gneqq0\quad \textrm{ in } \O\setminus \Gamma,
\end{cases}
\ee
where $q\in C^2(\O)$, $q\geq0$ in $\O$ and normalized as
\be\label{eq:ass-q}
 {\max_{\s\in \Gamma}q(\s)=1}.
\ee

We obtain the following result:

\begin{Theorem}\label{th:perturb}
Suppose that {$p\geq  \frac{N-k+2}{N-k-2}$} and that \eqref{eq:ass-q} holds.
 Then problem \eqref{eq:petplmi} does not have a solution.
\end{Theorem}
The above supercriticality assumption on $p$ is sharp as we will see in Section \ref{ss:exist-q} below. 

\begin{Remark}
\begin{itemize}
\item  It was observed in [\cite{BDT}, Remark 3] that  if $0<\max_{\G} q<1$ or $q\equiv1$ then \eqref{eq:petplmi} does not have a
solution when $$p\geq p^+:=1+\frac{2}{ \frac{N-k-2}{2}-\sqrt{\frac{(N-k-2)^2}{4}-c}  }\geq \frac{N-k+2}{N-k-2},$$ with $c=\frac{(N-k-2)^2}{4}\,\max_{\G} q$.

\item We should  mention that extremals for    weighted Hardy inequality was studied in  \cite{BM}, \cite{BMS} and \cite{FaMa}
when $\G$ is a submanifold of $\de\O$ and $k=1,\dots, N-1$. In these papers, the finiteness of
the integral
$ \displaystyle{\int_{\Gamma}\frac{1}{\sqrt{1-q(\s)}}\,d\s}$  
  was necessary and sufficient
to obtain the  existence of an eigenfunction  in some function space corresponding to some
``critical`` eigenvalue.\\
We belive that the argument in this paper and the results in \cite{FaMa} might be used to study
problem \eqref{eq:petplmi} but with $\G\subset \de\O$.
\end{itemize}
\end{Remark}

In the sequel, we denote by $\d(x):= \textrm{dist}(x\,,\,\Gamma)$. For $\b>0$,  we consider the interior of the tube around $\Gamma$ of
 radius $\b$ defined as  $\Gamma_\b:= \{ x\in \O\,:\,  \d(x)<\b \}$.
It is well known that if  $\b$ is positive small, the function $\d$ is smooth in $\Gamma_\b \setminus \Gamma$.
If $\b$ is small then  for all $x\in\Gamma_\b$, there exists a unique projection
 $\s(x)\in \G$  given by
\be\label{eq:repsigx}
\s(x)=x-\frac{1}{2}\n(\d^2)(x)=x-\d(x)\n\d(x).
\ee
In addition the function $\s$ is also smooth in $\Gamma_\b$, see for instance \cite{AS}.

%
%
From now on, we will consider $\b'$s for which  the projection function $\s$ is smooth.\\
Set
\be\label{eq:defo0}
\o_0(x)=\delta^{-\alpha (x)},
\ee
with
\be\label{eq:a-eq-aq}
\alpha(x)=\a_q(x)=\frac{N-k-2}{2}-\sqrt{\tilde \alpha (x)}
\ee
and where
$$
\tilde \alpha (x)=\left(\frac{N-k-2}{2}\right)^2\left(1- q(\sigma (x))+\delta (x)\right).
$$
Clearly $\a$ is well defined as soon as $q\leq 1$ on $\Gamma$.
Recall that $X_a(t)=|\log t|^a$, $t\in (0,1)$ and $a\in\R$. We define
$$
\o_a(x):=\o_0(x)\,X_a(\d(x)).
$$
We will need the following result which will be useful in the
 proof of Theorem \ref{th:perturb}.
 \begin{Lemma}\label{lem:sub-sup-per}
 Put $L_q:=-\D-\left(\frac{N-k-2}{2}\right)^2\,\d^{-2}\,q$. Then there exit $C,\b_0>0$
 depending only on $\Gamma$, $a$ and $\|q\|_{C^2(\ov{\O})}$  such that
\be\label{eq:om-a}
\left| L_q\,\o_a  -{2\,a\,\sqrt{\tilde\alpha}}{\delta^{-2}}\,X_{-1}\,\o_a
  +{a(a-1)}{\delta^{-2}} \,X_{-2}\,\o_a\right|\leq C|\log(\delta)|\,\delta^{-\frac32}\,\o_a,\quad\textrm{ in } \Gamma_{\b_0}.
\ee
 \end{Lemma}
\proof
We start by noticing that
\be\label{eq:Dwlog}
\Delta \o_0=\o_0\bigg( \Delta \log(\o_0)+|\nabla\log(\o_0)|^2 \bigg)
\ee
and that
\be\label{eq:Dlogw}
-\Delta \log(\o_0)=\Delta \alpha\log(\delta)+2\nabla\alpha\cdot \nabla
(\log(\delta))+\alpha\Delta \log(\delta).
\ee
We have
\be\label{eq:Dal}
-\Delta\alpha=\Delta\sqrt{\tilde \alpha}=\sqrt{\tilde \alpha}\,\left(\frac12 \Delta\log(\tilde \alpha)
+\frac14|\nabla \log(\tilde \alpha)|^2  \right).
\ee
By simple computations we get
$$
\sqrt{\tilde\alpha}  \nabla\log(\tilde\alpha)=\frac{\nabla\tilde\alpha}{\sqrt{\tilde\alpha}}=
\left(\frac{N-k-2}{2}\right)^2
\frac{-\nabla(q\circ\sigma)+\nabla\delta}{\sqrt{\tilde{\a}}}
$$
and
$$
\sqrt{\tilde\alpha}  \left|\Delta\log(\tilde\alpha)\right|\leq\frac{|\Delta\tilde\alpha|}{\sqrt{\tilde\alpha}}
+\frac{\left|  \nabla\tilde\alpha\right|^2}{\sqrt{\tilde\alpha}}.
$$
We deduce that
there exits a constant $\b_0>0$  depending only on $\G$ and $\|q\|_{C^2(\ov{\O})}$
such that
%
%
\be\label{eq:absD-alfa}
|\Delta\alpha|\le \frac{C}{\delta^{\frac32}}\quad\textrm{ in } \Gamma_{\b_0}.
\ee
Similar we have
\be\label{eq:nanlogd}
|\n \a\cdot\n\log\d|\le \frac{C}{\delta^{\frac32}}\quad\textrm{ in } \Gamma_{\b_0}.
\ee
Recall that (see for instance \cite{DD-H})
\be\label{eq:aDlogd}
\alpha\Delta\log(\delta)=\alpha\,\frac{N-k-2}{\delta^2}\,(1+O(\delta)).
\ee
Using \eqref{eq:Dal}, \eqref{eq:absD-alfa}, \eqref{eq:nanlogd} and \eqref{eq:aDlogd}
in the formula \eqref{eq:Dlogw}, we obtain the following estimate:
\be\label{eq:estDlogw}
\left| \Delta\log(\o_0)+ \alpha\,\frac{N-k-2}{\delta^2}\right| \leq   C\,\frac{|\log\d|}{\delta^{\frac32}} \quad\textrm{ in } \Gamma_{\b_0}.
\ee
We also have
$$
-\nabla(\log(\o_0))=\nabla(\alpha \log(\delta))=\alpha
\frac{\nabla\delta}{\delta}+\log(\delta)\nabla \alpha
$$
and thus
\be\label{eq:Nablogw}
\left||\nabla(\log(\o_0))|^2-\frac{\alpha^2}{\delta^2}\right|\leq
 C\,\frac{|\log\d|}{\delta^{\frac32}}\quad\textrm{ in } \Gamma_{\b_0}.
 \ee
By using \eqref{eq:estDlogw}, \eqref{eq:Nablogw} in the identity \eqref{eq:Dwlog}, we conclude that
\begin{eqnarray*}
\left| \frac{ \Delta \o_0 }{\o_0}+\alpha\,\frac{N-k-2}{\delta^2}  -\frac{\alpha^2}{\delta^2}  \right|
\leq C\,\frac{|\log\d|}{\delta^{\frac32}}\quad\textrm{ in } \Gamma_{\b_0}.
\end{eqnarray*}
We use the fact that $|q(x)-q(\s(x))|\leq C\d(x)$ to deduce  that
\begin{eqnarray*}
 { \Delta \o_0 }&=&
 \frac{(N-k-2)^2}{4}{\delta^{-2}}\,q(x)\,\o_0+O(|\log(\delta)|\,\delta^{-\frac32})\,\o_0 \quad\textrm{ in } \Gamma_{\b_0}.
\end{eqnarray*}
To conclude, we write
$$
\o_a(x):=\o_0(x)\left(-\log(\d(x))\right)^a
$$
and the proof of \eqref{eq:om-a} follows with some little computations. We skip the details.
\QED
\subsection{Proof of Theorem \ref{th:perturb}}
%
\textbf{Step I:} The following inequality holds:
\be\label{eq:coercv}
\int_{\Gamma_{\b_0}}|\n \vp|^2dx-\frac{(N-k-2)^2}{4}
\int_{\Gamma_{\b_0}} \d^{-2}q\vp^2dx\geq C \int_{\Gamma_{\b_0}} \d^{-2}X_{-2}\vp^2dx\quad
\ee
for any $ \vp\in C^\infty_c(\Gamma_{\b_0})$, with  $\b_0>0$ small depending only on   $K$ and $\|q\|_{C^2(\ov{\O})}$ and $C>0$ is a constant.\\
Indeed, observe that by \eqref{eq:om-a},
$$
-\frac{\D\o_\frac{1}{2}}{\o_{\frac{1}{2}}}-\frac{(N-k-2)^2}{4}\delta^{-2}\,q\geq
\frac{1}{4}{\delta^{-2}} \,X_{-2}-C|\log(\delta)|\,\delta^{-\frac32}\quad\textrm{ in }\Gamma_\b\setminus \Gamma.
$$
Hence, there exist  $\b_0>0$ small and  a constant $C>0$ such that
\be\label{eq:w12geq}
-{\D\o_\frac{1}{2}}-\frac{(N-k-2)^2}{4}\delta^{-2}\,q\,{\o_{\frac{1}{2}}}
-C {\delta^{-2}} \,X_{-2}\,{\o_{\frac{1}{2}}}\geq0\quad\textrm{ in }\Gamma_{\b_0}\setminus \Gamma.
\ee
Since $\o_{\frac12}\in L^1(\Gamma_{\b_0})$, the inequality \eqref{eq:w12geq} holds in $\calD'(\Gamma_{\b_0})$ thus
 by Lemma \ref{lem:AP}, \eqref{eq:coercv} follows.\\
\bigskip
\noindent
\textbf{Step II:} Set $\th_a:=\o_0+\o_a$, with $a<-1/2$. There exist  positive constants $C$ and  $\b_0$ depending only on $a,$ $\G$ and $\|q\|_{C^2(\ov{\O})}$
such that  
\be\label{eq:normtheta-H1}
\|\th_a\|^2
_{H^1(\Gamma_{\b_0})}\leq C \displaystyle{\int_{\Gamma}\frac{1}{\sqrt{1-q(\s)}}\,d\s}.
,
\ee

\bigskip
\noindent
First of all it is easy to see that,
 since $X_a\leq1$ for $a$ negative, we can estimate
\be\label{eq:estnth}
|\n \th_a|^2
\leq C \d^{-2\a-2}\quad\textrm{ in } \Gamma_{\b_0}.
\ee
Following \cite{DD}, there exits a family of disjoint open sets $W_i$, $i=1,\dots,m_0$ of $\Gamma$ such that
$$
\Gamma=\bigcup_{i=1}^{m_0}\ov{W_i},\qquad |\ov{W_i}\cap \ov{W_j} |=0,\quad i\neq j.
$$
Moreover by \eqref{eq:estnth},
\be\label{eq:polarcordCylin}
  \|\th_a\|^2
_{H^1(\Gamma_{\b_0})}\leq C  \int_{\Gamma_{\b_0}}\d^{-2\a-2}=C\sum_{i=1}^{m_0} \int_{W_i\times B^{N-k}_{\b_0}}\d^{-2\a-2}\,(1+O_i(\d))\,d\d\,d\s,
\ee
where  $B^{N-k}_\b$ is the ball of $\R^{N-k}$ with radius $\b$. Therefore, we have
\begin{eqnarray*}
\|\th_a\|^2
_{H^1(\Gamma_{\b_0})} &\leq& C  \sum_{i=1}^{m_0} \int_{W_i}\int_{\S^{N-k-1}}\int_0^{\b_0}
\d^{-1}\d^{(N-k-2)\sqrt{1-q(\s)+\d}}\,d\d d\s\\
&\leq &C\sum_{i=1}^{m_0} \int_{W_i}\int_{\S^{N-k-1}}\int_0^{\b_0}
\d^{-1}\d^{(N-k-2)\sqrt{1-q(\s)}}\,d\d d\s\\
&\leq &C\sum_{i=1}^{m_0} \int_{W_i}\frac{1}{\sqrt{1-q(\s)}}\,d\s=C\int_{\Gamma}\frac{1}{\sqrt{1-q(\s)}}\,d\s.
\end{eqnarray*}
This ends the proof of this step.\\
\bigskip
\noindent
\textbf{Step III}: Let $u$ satisfies \eqref{eq:petplmi} and $\th_a=\o_0+\o_a$, for $a<-1/2$.
For any $\b>0$ small, there exists a constant $C>0$ such that
\be\label{eq:lestu}
u\geq C \th_a\quad\textrm{ in } \Gamma_\b .
\ee
Indeed, define $q_n(x):=q(x)-\frac{1}{n}$ with $n\in \N^*$ and
 we put $\th_{a,n}=\d^{-\a_{q_n}}+\d^{-\a_{q_n}}X_a(\d)$.
Recalling  \eqref{eq:a-eq-aq}, by \eqref{eq:om-a}  there exit  constants $\b_0, C>0$ (independent on $n$) such that
$$
L_{q_n}\th_{a,n}\leq
 - \frac{3}{4}{\delta^{-2}} \,|\log\d|^{-2+a}\,\d^{-\a_{q_n}}+ C|\log(\delta)|\,\delta^{-\frac32}\,\d^{-\a_{q_n}}\quad \textrm{ in } \Gamma_\b,
$$
for any  $\b\in(0,\b_0)$.
Therefore for   all $\b>0$ small   we obtain
\be\label{eq:thak-suh}
-\D\th_{a,n}-\frac{(N-k-2)^2}{4}\delta^{-2}\,q_n(x)\,\th_{a,n}\leq0\quad \textrm{ in }\Gamma_\b\quad \forall  n\geq1.
\ee
By [\cite{BDT}, Lemma 1], $u\in L^p_{loc}(\O)$. In addition, it is  nonnegative and non-trivial in   $\O$ and satisfies
\be\label{eq:urmsing}
-\D u- \frac{(N-k-2)^2}{4}\d^{-2}q(x)\, u\geq u^p\quad\textrm{ in }\calD'(\O).
\ee
Hence by the maximum principle, $u>0$ in $\O$. For $\b>0$  small (independent on $n$), by \eqref{eq:coercv} 
we can pick $v_n\in H^1_0(\Gamma_\b)$   solution to
\be\label{eq:vq}
-\D v_n- \frac{(N-k-2)^2}{4}\d^{-2}q_n(x)\, v_n= \min(u^p,1)\quad\textrm{ in }\Gamma_\b.
\ee
By Lemma \ref{lem:min-sol} the sequence $(v_n)_n$ is monotone increasing and converging pointwise 
to $v\in H(\G_\b)$ solution to $-\D v- \frac{(N-k-2)^2}{4}\d^{-2}q(x)\, v= \min(u^p,1)$.
By Lemma \ref{lem:min-sol}
 we have that $u\geq v\geq v_n>0$ in $\Gamma_\b$ for any $n\geq1$. 
By elliptic regularity theory $v_n$ is continuous in $ \Gamma_\b\setminus \Gamma$.
 We choose  $M_n>0$ such that
\be\label{eq:Mk}
M_n\sup_{\de\Gamma_{\frac{\b}{2}}} \th_{a}=\inf_{\de\Gamma_{\frac{\b}{2}}} v_{n}.
\ee
Clearly, we have $M_n\,\th_{a,n}\leq v_{n}$ on $ \de\Gamma_{\frac{\b}{2}}$.
It follows form  \eqref{eq:normtheta-H1} that $\left(M_n\,\th_{a,n}- v_{n}\right)^+\in H^1_0\left( \Gamma_{\frac{\b}{2}}\right)$.
On the other hand by \eqref{eq:thak-suh} and \eqref{eq:vq},
$$
-\D \left( M_n\,\th_{a,n}- v_{n}\right)- \frac{(N-k-2)^2}{4}\d^{-2}q(x) \left( M_n\,\th_{a,n}- v_{n}\right)\leq0\quad\textrm{ in }\Gamma_{\frac{\b}{2}}.
$$
Multiplying this inequality by $\left(M_n\,\th_{a,n}- v_{n}\right)^+ $ and integrating by
 parts yields $M_n\,\th_{a,n}\leq v_{n}$ on $ \Gamma_{\frac{\b}{2}}$ by \eqref{eq:coercv}. Since $v_{n}$ is
 monotone increasing to $v$,  by the choice of $M_n$ in \eqref{eq:Mk}, there exists an integer $n_0\geq1$ such that $M_{n_0}\,\th_{a,n}\leq v_{n}$ for all $n\geq n_0$.
 Passing to the limit, we get \eqref{eq:lestu}.\\
\bigskip
\noindent
\textbf{Step IV:}
There is no  $u$ satisfying \eqref{eq:petplmi} with $p\geq\frac{N-k+2}{N-k-2}$.\\
\bigskip
\noindent
By using  \eqref{eq:lestu}  we have  that
$$
 u^{p-1}\geq C \th_a^{p-1}\geq C\o_0^{p-1}\geq C \d^{-2+2\sqrt{1-q\circ\s}}\quad \textrm{ in } \G_\b,
$$
for some $C>0$ and  provided $\b$ is small.
 This together with   \eqref{eq:urmsing} give
\be\label{eq:lin-eq-q}
-\D u -\left(q+C_0\d^{2\sqrt{1-q\circ\s}}\right)\frac{(N-k-2)^2}{4}\d^{-2}\,u\geq0 \qquad  \textrm{ in } \calD'(\G_\b),
\ee
for some $C_0>0$. 
By Lemma \ref{lem:AP} we have, $\forall \vp\in C^\infty_c(\G_\b)$
\be\label{eq:aplAP-q}
\frac{(N-k-2)^2}{4}\leq \frac{\displaystyle \int_{\G_\b}|\n \vp|^2\,dx
}{\displaystyle
\int_{\G_\b}\left(q+ C_0\d^{2\sqrt{1-q(\s)}}  \right)\,\d^{-2}\,\vp^2\,dx}.
\ee
Our aim is to construct appropriate  test  functions in \eqref{eq:aplAP-q} supported in a neighborhood of the maximum point of $q$ on $\G$ in order to get 
a contradiction.\\
 
By \eqref{eq:ass-q}, we can let $\sigma_0\in\G$ be such that
\be\label{eq:maxq-eq1}
q(\sigma_0)=\max_{\s\in \G}q(\s)=1.
\ee

For $y\in\R^N$, we write $y=(\tilde{y},\bar{y}) \in \R^{N-k}\times \R^k$ with 
$\tilde{y}=(y^1,\dots, y^{N-k})$ and $\bar{y}=(y^{N-k+1},\dots,y^N)$.
 Consider  $f:\R^k\to \G $   a normal parameterization of  a neighborhood of $\s_0$ with $f(0)=\s_0$.
In a    neighborhood of $\s_0$, we consider
 $\mathcal{N}_i$, $i=1,\dots,N-k$ an orthonormal frame filed on the normal bundle of $\G$. We can therefore define a parameterization 
of a neighborhood, in $\R^N$, of $\s_0$ by the mapping $Y:B_r(0)\to \G_\b$   as 
$$
  y \mapsto Y(y)=f(\bar{y})+\sum_{i=1}^{N-k}y^i\mathcal{N}_i(f(\bar{y}))\in \G_\b,
$$
 for some $r>0$ small.
By identification using \eqref{eq:repsigx}, we get  for some $r>0$ small
\be\label{eq:ident-q}
\d( Y(y))=|\tilde{y}|,\qquad \s(Y(y))=f(\bar{y}) \qquad \forall y\in B_r(0).
\ee
Denoting by $g$    the metric induced by $Y$ with component $g_{ij}(y)=\la \de_i Y(y),\de_j Y(y)\ra$,
it is not difficult to verify that  for all $ y\in B_r(0)$
\be\label{eq:gij}
g_{ij}(y)=\d_{ij}+O(|y|)\qquad \textrm{ for } i,j=1,\dots,N.
\ee

Next we let $w\in C^\infty_c(\R^{N-k}\setminus\{0\}\times \R^k)$.
We choose
$\e_0>0$ small such that, for all $\e\in(0,\e_0) $, we have  
$$
\e\,\textrm{Supp} w\subset B_r(0).
$$
We define the following test function
$$
\vp_\e(x)=\e^{\frac{2-N}{2}}w\left(\e^{-1}Y^{-1}(x)\right),
\quad x\in Y(\e\,\textrm{Supp} w ).
$$
Clearly, for every $\e\in(0,\e_0)$, we have that $\vp_\e\in
C^\infty_c(\G_\b)$  and thus by \eqref{eq:aplAP-q},  we have (summations over repeated indices is understood)
\begin{align*}
\frac{(N-k-2)^2}{4}&\leq\frac{\displaystyle \int_{\G_\b}|\n \vp_\e|^2\,dx
}{\displaystyle
\int_{\G_\b}\left(q+ C_0\d^{2\sqrt{1-q(\s)}} \right)\,\d^{-2}\,\vp_\e^2\,dx}\hspace{10cm}\\ 
&=\frac{ \e^{2-N}\displaystyle
\int_{\R^N}\e^{-2}(g^\e)^{ij}\partial_i
w\partial_j w\,\sqrt{|g^\e}| \,dy
}{ \,\displaystyle
\e^{2-N}\int_{\R^N}\left(q(Y(\e y))+ C_0|\e \tilde{y}|^{2\sqrt{1-q( f(\e \bar{y}))}}  \right)\,|\e\tilde y|^{-2}\,w^2\,\sqrt{|g^\e|} \,d y}  ,\\
&=\frac{  \displaystyle
\int_{\R^N} (g^\e)^{ij}\partial_i
w\partial_j w\,\sqrt{|g^\e}| \,dy
}{ \,\displaystyle
\int_{\R^N} \,\left(q(Y(\e y))+ C_0|\e \tilde{y}|^{2\sqrt{1-q( f(\e \bar{y}))}}  \right)\,|\tilde y|^{-2}\,w^2\,\sqrt{|g^\e|} \,d y}  ,\\
\end{align*}
where $g^\e$ is the metric with component  $g^\e_{ij}(y)= g_{ij}(\e y) $
with  $(g^\e)^{ij}(y)$ denotes the component of the inverse matrix of $g^\e$  and   $|g^\e|$ stands for the determinant  of $g_\e$.

Observe that  the scaled metric $g^\e$ expands a $g^\e=Id+O(\e)$ on the support of $w$ by \eqref{eq:gij}. In addition since $q$ is
of class $C^1$, decreasing $\e_0$ if necessary, there exits $c>0$ such that 
$$
1-q( f(\e \bar{y}))\leq c \e \qquad \forall \bar{y}\in \textrm{Supp}w\cap \R^k,\quad\forall \e\in(0,\e_0),
$$
 by \eqref{eq:ident-q}. From this we deduce that 
$$
|\e \tilde{y}|^{2\sqrt{1-q( f(\e \bar{y}))}}\to 1 \quad \textrm{  as } \e\to0,
$$
uniformly in ${y}\in \textrm{Supp}w $.
We then  have  from the dominated convergence theorem  and  using \eqref{eq:maxq-eq1}
together with  \eqref{eq:ident-q} 
\begin{eqnarray*}
\frac{(N-k-2)^2}{4} &\le&  \frac{1}{1+C_0} \frac{\displaystyle
\int_{\R^N}|\nabla w|^2\,d y }{\displaystyle
\int_{\R^N} |\tilde y|^{-2}\,w^2\,d y}\qquad \forall w\in C^\infty_c(\R^{N-k}\setminus\{0\}\times \R^k) .
\end{eqnarray*}
This is in contradiction with the well  know fact that 
$$
 \inf_{ w\in C^\infty_c(\R^{N-k}\setminus\{0\}\times \R^k) }
 \frac{\displaystyle \int_{\R^N}|\nabla w|^2\,d y }{\displaystyle
\int_{\R^N} |\tilde y|^{-2}\,w^2\,d y}=  \inf_{ w\in C^\infty_c(\R^{N}) }
 \frac{\displaystyle \int_{\R^N}|\nabla w|^2\,d y }{\displaystyle
\int_{\R^N} |\tilde y|^{-2}\,w^2\,d y}= \frac{(N-k-2)^2}{4}
$$
because $N-k>2$, see for instance [\cite{Maz},    Section 2.1.6] and [\cite{Musina}, Lemma 1.1].
\QED
%
\subsection{Existence}\label{ss:exist-q}
\begin{Proposition}\label{prop:exist-delta}
Let $1\leq  p<\frac{N-k+2}{N-k-2}$.
Then if $\b$ is small, there exists $u\in L^p(\G_\b)$ satisfying
\be\label{eq:petplmi-super}
\begin{cases}
-\D u - \frac{(N-k-2)^2}{4}\,{\d^{-2}}\, {q}\,u\geq u^p
\quad\textrm{ in } \G_\b\setminus \Gamma, \\
u>0\quad \textrm{ in } \G_\b.
\end{cases}
\ee
\end{Proposition}
\proof
Set 
$$
u=\o_0-\o_{-1}=\o_0(1-X_{-1}(\d)).
$$
Then by Lemma \ref{lem:sub-sup-per}
there exits $ C>0$  such that
$$
L_{q}u \geq
2{\delta^{-2}} \,X_{-3}(\d)\,\d^{-\a}- CX_{1}(\d)\,\delta^{-\frac32}\,\d^{-\a}
\quad \textrm{ in } \Gamma_\b\setminus\G.
$$
Hence, provided $\b$ is small, we have $u>0$ and
$$
-\D u - \frac{(N-k-2)^2}{4}\,{\d^{-2}}\, {q}\,u\geq 
 \d^{-2}X_{-5}(\d)u\quad \textrm{ in } \G_\b\setminus\G.
$$
We thus want 
$$
\d^{-2}X_{-5}(\d)u\geq u^p\quad \textrm{ in } \G_\b\setminus\G.
$$
Or equivalently
$$
\d^{-2}X_{-5}(\d)\d^{-\a}(1-X_{-1}(\d))\geq \d^{-p\a}(1-X_{-1}(\d))^p
\quad \textrm{ in } \G_\b\setminus\G.
$$
That is 
\be\label{eq:needd}
\d^{-2+(p-1)\a}X_{-5}(\d)(1-X_{-1}(\d))^{1-p}\geq 1 \quad \textrm{ in } \G_\b\setminus\G.
\ee
We observe that for $1\leq  p< \frac{N-k+2}{N-k-2}$ we have for every $x\in \G_\b\setminus\G$
$$
-2+(p-1)\a(x)\leq -2 -(p-1)\frac{N-k-2}{2} 
<0.
$$
This implies that if $\b$ is small enough,  \eqref{eq:needd} holds so that $u$
satisfies \eqref{eq:petplmi-super}. The fact that $u\in L^p(\G_\b)$ is easy to check, we skip the details.
\QED
\textbf{ Acknowledgments} \\
This work is supported by the Alexander-von-Humboldt Foundation.
The author would like to thank Professor Roberta Musina and Professor
Tobias Weth for their useful comments and suggestions. We gratefully  thank the referee for 
carefully reading this manuscript and who's comments have 
helped to improve considerably the former versions  and also for taking  the author's attention to \cite{CFKS}.

 \label{References}


\begin{thebibliography}
\footnotesize
%
\bibitem{AS} L. Ambrosio and H. M. Soner, Level set approach to mean curvature flow in arbitrary
codimension. J. Diff. Geometry, 43, (1996) 693-737.

\bibitem{AS-CPDE} S. Armstrong, B. Sirakov,   Nonexistence of positive supersolutions of elliptic equations via the maximum principle.
Comm. PDE, 36, 11 (2011) 2011-2047.
\bibitem{AS-RIMS} S. Armstrong, B. Sirakov, A new approach to Liouville theorems for elliptic inequalities. RIMS proceeding, to appear.
%
%
\bibitem{BC} H. Brezis and X. Cabr\'e. Some simple nonlinear PDEs without solutions. Bull. UMI 1 (1998),
223-262.
%
\bibitem{BDT}H. Brezis, L. Dupaigne, A. Tesei,
On a semilinear elliptic equation with inverse-square potential.
Selecta Math. (N.S.) 11 (2005), no. 1, 1-7.
%
\bibitem{BM} H. Brezis and M. Marcus, Hardy's inequalities revisited.
 Dedicated to Ennio De Giorgi.
  Ann. Scuola Norm. Sup. Pisa Cl. Sci. (4)  25  (1997),  no. 1-2, 217-237.
%
 \bibitem{BMS} H. Brezis, M. Marcus and I. Shafrir,
   Extermal functions for Hardy's inequality with weight, J. Funct. Anal. 171 (2000), 177-191.
%
\bibitem{CFKS} H. L. Cycon, R. G. Froese, W. Kirsch, and B. Simon, Schr\"{o}dinger operators with application
to quantum mechanics and global geometry, Berlin, 1987.
%
\bibitem{DD-H}J. D\'{a}vila, L. Dupaigne,  Hardy-type inequalities.
 J. Eur. Math. Soc. (JEMS) 6 (2004), no. 3, 335-365.
%
\bibitem{DD}J. D\'{a}vila, L.  Dupaigne, Comparison results for PDEs with a singular potential.
 Proc. Roy. Soc. Edinburgh Sect. A  133  (2003),  no. 1, 61-83.
%
\bibitem{Docarmo} M. P.  do Carmo, Riemannian geometry, Mathematics: Theory and Applications. Birkh\"{a}user Boston, Inc., Bostonn, MA, 1992.
%
%
\bibitem{D}L. Dupaigne, Semilinear elliptic PDEs with a singular potential, J. Anal. Math., 86
(2002), 359-398.
%
%
\bibitem{DN} L. Dupaigne and G. Nedev. Semilinear elliptic PDEs with a singular potential. Adv. Dif-
ferential Equations 7 (2002), 973-1002.
%
\bibitem{FaMa}M. M. Fall, F. Mahmoudi,  Weighted Hardy inequality with higher dimensional singularity on the boundary.
http://arxiv.org/abs/1202.0033.
%
\bibitem{mmf}M. M. Fall, On the Hardy Poincar\'e inequality with boundary singularities.
Commun. Contemp. Math., Vol. 14, No. 3 (2012) 1250019.
%
\bibitem{mmf-p} M. M. Fall, A note on Hardy's inequalities with
boundary singularities. Nonlinear Anal. 75 (2012) no. 2, 951--963.
%
%
\bibitem{FaMu1}M. M. Fall, R. Musina,  Hardy-Poincar\'e inequality with boundary singularities.
Proc. Roy. Soc. Edinburgh. A 142 (2012) 1--18.
%
\bibitem{FaMu-ne} M. M. Fall, R. Musina,  Sharp nonexistence results for a linear elliptic
 inequality involving Hardy and Leray potentials.
J. Inequal.  Appl. 2011 (2011). doi:10.1155/2011/917201.
%
%
\bibitem{MMF-PJM} M. M. Fall, Area-minimizing regions with small volume
in Riemannian manifolds with boundary. Pacific J. Math. 244 (2010),
no. 2, 235-260.
%

%
\bibitem{KLS-Trans}V. Kondratiev, V. Liskevich,  Z.  Sobol, Positive supersolutions to semi-linear second-order non-divergence type elliptic equations in exterior domains.
Trans. Amer. Math. Soc. Volume 361 (2009), 697--713.
%
\bibitem{KLS}V. Kondratiev, V. Liskevich,  Z.  Sobol,
Positive solutions to semi-linear and quasi-linear second-order elliptic equations on unbounded domains. Handbook of  Differential  Equation.
Stationary Partial Differential Equations,
volume 6 (Edited by M. Chipot) 2008 Elsevier, pp.177-268.
%
\bibitem{LLM} V. Liskevich, S. Lyakhova and V. Moroz, Positive solutions to singular semilinear elliptic
equations with critical potential on cone-like domains, Adv. Differential Equations 11 (2006), pp. 361-398.
%
%
\bibitem{LLM-Proce} V. Liskevich, S. Lyakhova and V. Moroz, Positive solutions to semilinear elliptic equations with critical lower order terms.
EQUADIFF 2003, 549–554, World Sci. Publ., Hackensack, NJ, 2005.
%
\bibitem{KLS-JDE}V. Kondratiev, V. Liskevich  Z.  Sobol, Second-order semilinear elliptic inequalities in exterior domains.
J. of Differential Equations, vol. 187 (2003), 429-455.
%
%
\bibitem{Maz} V. G. Maz'ja, Sobolev Spaces, Springer-Verlag, Berlin, 1980.
%
\bibitem{Musina} R. Musina, Ground state solutions of a critical problem involving
cylindrical weights. Nonlinear Analysis 68 (2008) 3972-3986.
%
\bibitem{PoTe}S. I. Pohozaev and A. Tesei. Nonexistence of local solutions to semilinear partial differential
inequalities. Ann. Inst. H. Poincar\'e Anal. Non Lin\'eaire 21 (2004), pp. 487-502.
%
 \bibitem{PT} Y. Pinchover\ and\ K. Tintarev, Existence of minimizers for Schr\"odinger operators
under domain perturbations with application to Hardy's inequality, Indiana Univ. Math. J. {\bf 54} (2005), no. 4, 1061--1074.
%
\bibitem{T} S. Terracini, On positive entire solutions to a class of equations with a singular coefficient and critical exponent.
 Adv. Differential Equations  1  (1996),  no. 2, 241--264.
%
\end{thebibliography}
\end{document}